\documentclass[12pt]{amsart}
\usepackage{amsfonts}
\usepackage{amssymb, amstext, amscd, amsmath}
\usepackage{amsthm}
\usepackage{times}
\usepackage{txfonts}
\usepackage{indentfirst}
\usepackage[all]{xy}
\usepackage{subfigure}
\usepackage{graphicx}
\newtheorem{thm}{Theorem}[section]

\newtheorem{lem}[thm]{Lemma}
\newtheorem{conj}[thm]{Conjecture}

\theoremstyle{definition}

\newtheorem{dfn}[thm]{Definition}

\numberwithin{equation}{section}
\numberwithin{equation}{section}
\def\square{\vbox{
      \hrule height 0.4pt
      \hbox{\vrule width 0.4pt height 5.5pt \kern 5.5pt \vrule width 0.4pt}
      \hrule height 0.4pt}}

\def\ch\mathrm{c h}

\begin{document}

\title[The Slope Conjecture for a Family of Montesinos Knots]
{The Slope Conjecture for a Family of Montesinos Knots}

\author{Xudong Leng$^{*}$}
\address{School of Mathematical Sciences, Dalian University of
Technology, Dalian 116024, P. R. China} \email{xudleng@163.com}
\thanks{$^{*}$ Corresponding author}
\author{Zhiqing Yang$^{\dag}$}
\address{School of Mathematical Sciences, Dalian University of
Technology, Dalian 116024, P. R. China} \email{yangzhq@dlut.edu.cn}
\thanks{$^{\dag}$ Supported by the NFSC (No. 11271058)}
\author{Xinmin Liu$^{\ddag}$}
\address{School of Mathematical Sciences, Dalian University of
Technology, Dalian 116024, P. R. China }
\email{ximinliu@dlut.edu.cn}
\thanks{$^{\ddag}$ Supported by the NFSC (No. 10371076 and 11431009) }



\begin{abstract}
The Slope Conjecture relates the degree of the colored Jones polynomial to the boundary slopes of a knot. We verify the Slope Conjecture and the Strong Slope Conjecture for Montesinos knots $M(\frac{1}{r},\frac{1}{s-\frac{1}{u}},\frac{1}{t} )$ with $r,u,t$ odd, $s$ even and $u\leq-1$, $r<-1<1<s,t$.
\end{abstract}
\maketitle


\section{Introduction}
Soon after V. Jones discovered the famous Jones polynomial~\cite{J}, E. Witten found an intrinsic explanation~\cite{W} through TQFT approach, which led to the colored Jones polynomial. As a generalization of Jones polynomial, colored Jones polynomial reveals many deep connections between quantum algebra and three-dimensional topology. For example, the Volume Conjecture, which relates the asymptotic behavior of the colored Jones polynomial of a knot to the hyperbolic volume of its complement. Another connection proposed by S. Garoufalidis~\cite{Gar11b} named Slope Conjecture, predicts that the growth of maximal degree of colored Jones polynomial of a knot determines some boundary slopes of the knot complement. So far, the Slope Conjecture has been verified for knots with up to 10 crossings~\cite{Gar11b}, adequate knots~\cite{FKP11}, 2-fusion knots~\cite{GvdV} and some pretzel knots~\cite{LV}. In~\cite{MT15}, K. Motegi and T. Takata prove that the conjecture is closed under taking connected sums and is true for graph knots. In~\cite{KT15}, E. Kalfagianni and A. T. Tran show the conjecture is closed under taking the $(p,q)$-cable with certain conditions on the colored Jones polynomial and formulate the Strong Slope Conjecture, see Conjecture 2.2b.

Enlightened by Lee and van der Veen~\cite{LV}, in this article we prove the Slope Conjecture and the Strong Slope Conjecture for a family of Montesinos knots, $M(\frac{1}{r},\frac{1}{s-\frac{1}{u}},\frac{1}{t} )$ with $r,u,t$ odd, $s$ even and  $u\leq -1$, $r<-1<1<s,t$ (See Section 2 and Figure 1). Particularly, when $u=-1$, $M(\frac{1}{r},\frac{1}{s+1},\frac{1}{t} )$ are just the pretzel knots in~\cite{LV}, and our results coincide with that of~\cite{LV} in this case. The strategy of the proof is to compare the maximal degree of colored Jones polynomial and certain boundary slopes of the corresponding knot. To compute colored Jones polynomial, we use the notion of knotted trivalent graphs~\cite{LV,vdV09,Thu02}, which is a convenient version of skein method~\cite{MV} for Montesinos knots, and to compute boundary slope we apply Hather and Oertel's edgepath system~\cite{HO89}.

\section{The Slope Conjecture}
Let $K$ denote a knot in $S^3$ and $N(K)$ denote its tubular neighbourhood. A surface $S$ properly embedded in the knot exterior $E(K)=S^3-N(K)$ is called \textit{essential} if it is incompressible, $\partial$-incompressible, and non $\partial$-parallel. A fraction $\frac{p}{q}\in \mathbb{Q}\bigcup \{\infty\}$ is a \textit{boundary slope} of $K$ if $pm+ql$ represents the homology class of $\partial S$ in the torus $\partial N(K)$, where $m$ and $l$ are the canonical meridian and longitude basis of $H_1(\partial N(K))$. The \textit{number of sheets} of $S$, denoted by $\sharp S$, is the minimal number of points at which the meridional circle of $\partial N(K)$ and $\partial S$ intersect.

For colored Jones polynomial, we use the convention of~\cite{LV}, where the unnormalized $\textit{n-colored Jones polynomial}$ is denoted by $J_{K}(n;v)$, see Section 3. Its value on the unknot is $[n]=\frac{v^{2n}-v^{-2n}}{v^2-v^{-2}}$ and the variable $v$ satisfies $v= A^{-1}$, where A is the A-variable of the Kauffman bracket. The maximal degree of $J_{K}(n)$ in $v$ is denoted by $d_{+}J_{K}(n)$.

A fundamental result due to S. Garoufalidis and T. Q. T. Le~\cite{GL05} states that colored Jones polynomial is \textit{q}-holonomic. Furthermore, the degree of a colored Jones polynomial is a quadratic quasi-polynomial~\cite{Gar11a}, which can be formulated as follows.

\begin{thm}~\cite{Gar11a}
For any knot $K$, there exist integer $p_{K}\in \mathbb{N}$ and quadratic polynomials $Q_{K,1}\ldots Q_{K,j} \in \mathbb{Q}[x]$ such that $d_{+} J_{K}(n)=Q_{K,j}(n)$ if $n=j\ (mod\ p_{K})$ for n sufficiently large.
\end{thm}

Now we can state the Slope Conjecture and the Strong Slope Conjecture:
\begin{conj}
In the context of the above theorem, set $Q_{K,j}=a_{j}x^2+ 2b_{j}x+ c_j$, then for each $j$ there exists an essential surface $S_j \subset S^3-K$, such that:

a.(Slope Conjecture~\cite{Gar11b}) $a_j$ is a boundary slope of $S_j$,

b.(Strong Slope Conjecture~\cite{KT15}) $ b_j=\frac{\chi(S_j)}{\sharp S_j}$, where $\chi(S_j)$ is the Euler characteristic of $S_j$.

\end{conj}

\begin{figure}[!ht]
\centering

\includegraphics{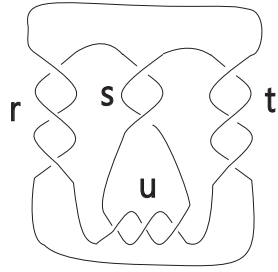}

\caption{The Montesinos Knot $M(\frac{1}{r},\frac{1}{s-\frac{1}{u}},\frac{1}{t} )$ with $r=u=-3, s=2, t=3$}
\end{figure}

A Montesinos knot is defined as a knot obtained by putting rational tangles together in a circle (See figure 1). A Montesinos knot obtained from rational tangles $R_{1}$, $R_{2}$, \ldots $R_{N}$ is denoted by $ M(R_{1}, R_{2},..., R_{N})$. Properties about Montesinos knots are omitted here, for details see~\cite{BZ}. It is known that a Montesinos knot is always semi-adequate~\cite{LT}, and the Slope Conjecture has been proved for adequate knots~\cite{FKP11}. So we focus on a family of A-adequate and non-B adequate knots $M(\frac{1}{r},\frac{1}{s-\frac{1}{u}},\frac{1}{t} ) $ with $r,u,t$ odd, $s$ even and $u\leq-1$, $r<-1<1<s,t$ and the maximal degrees of their colored Jones polynomials. The following is our main theorem.

\begin{thm}
The Slope Conjecture and the Strong Slope Conjecture are true for the Montesinos knots $M(\frac{1}{r},\frac{1}{s-\frac{1}{u}},\frac{1}{t} ) $, where $r,u,t$ are odd and $s$ is even and $u\leq -1$, $r<-1<1<s,t$.

\end{thm}

The above theorem is proved directly from the following two theorems. The first one is about the degree of the colored Jones polynomial and the second is about essential surfaces.

\begin{thm}
Let $K=M(\frac{1}{r},\frac{1}{s-\frac{1}{u}},\frac{1}{t} )$, where $r,u,t$ are odd, $s$ is even and $u\leq-1$, $r<-1<1<s,t$,
, set $A=-\frac{r+s+1}{2}, B=-(r+1), C=-\frac{r+t}{2}, \Delta=4AC-B^2 $.  \\
(1)If $ A\geq0 $ or $ C\geq0 $, or $ A,C<0 \ and\ \Delta<0 $, then
\[d_{+} J_{K}(n)= Q_{K,j}=[\frac{2(t-1)^{2}}{s+t-1}-2(r+t)]n^{2}+2(r+u+3)n+c_{j},\]
where $c_j$ is defined as following. Let $0\leq j < \frac{s+t-1}{2}$ such that  $n= j\ mod\  \frac{s+t-1}{2}$ , and set $v_j$ to be the odd number nearest to $\frac{2(t-1)j}{s+t-1} $. Then $c_j=-\frac{s+t-1}{2}\beta^{2}_{j}-(s+t-1)\beta_j-2(u+2)$, $\beta_j=v_j-1-\frac{2(t-1)}{s+t-1}j$. $p_K=\frac{s+t-1}{2}$ is a period of $d_{+} J_{K}(n)$ but may not be the least one. \\
(2)If $A,C< 0 \ and \ \Delta\geq0$, then
\[d_{+} J_{K}(n)= 2u(n-1).\]

\end{thm}

\begin{thm}
Under the same assumptions as the previous theorem, \\
(1)When $ A\geq0 $ or $ C\geq0 $, or $ A,C<0 \ and\ \Delta<0 $, there exists an essential surface $S$ with boundary slope $\frac{2(t-1)^{2}}{s+t-1}-2(r+t)$, and $\frac{\chi_{(S)}}{\sharp S}=r+u+3$.\\
(2)When $A,C< 0 \ and \ \Delta\geq0$, there exist an essential surface $S_0$ with boundary slope $0$, and $\frac{\chi_{(S_0)}}{\sharp S_0}= u$.
\end{thm}

\section{The Colored Jones Polynomial}
To compute the colored Jones polynomial of Montesinos knots, we use the notion of \textit{knotted trivalent graphs} (KTG) introduced in~\cite{LV} (See also~\cite{vdV09,Thu02}.), which is a natural generalization of knots and links.

\begin{figure}[!ht]
\centering

\includegraphics{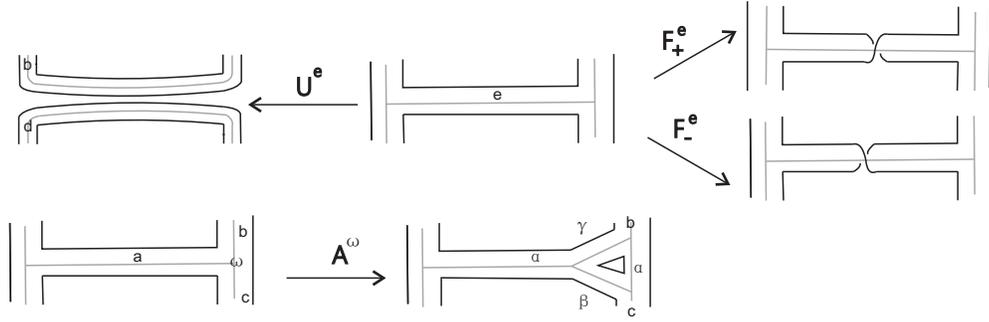}

\caption{Operations on KTG: framing change $F$ and unzip $U$ applied to an edge $e$, triangle move $A^{\omega}$ applied to a vertex $\omega$.}
\end{figure}

\begin{dfn}
~\cite{LV}

(1)A \textit{framed graph} is a one dimensional simplicial complex $\Gamma$ together with an embedding $\Gamma \rightarrow \Sigma$ of $\Gamma$ into a  surface with boundary of $\Sigma$ as a spine.

(2)A \textit{coloring} of $\Gamma$ is a map $\sigma: E(\Gamma)\rightarrow\mathbb{N}$, where $E(\Gamma)$ is the set of edges of $\Gamma$.

(3)A \textit{knotted trivalent graph} (KTG) is a trivalent framed graph embedded as a surface into $\mathbb{R}^3$, considered up to isotopy.
\end{dfn}
The advantage of KTGs over knots or links is that they support many operations. There are three types of operations we will need in this paper, the \textit{framing change} denoted by $F^{e}_{\pm}$, the \textit{unzip} denoted by $U^e$, and the \textit{triangle move} denoted by $A^{\omega}$, as illustrated in Figure 2.

The important thing is that these three types of operations are sufficient to produce any KTG from the $\Theta$ graph.

\begin{thm}
~\cite{vdV09,Thu02} Any KTG can be generated from $\Theta$ by repeatedly applying the operations $F_{\pm}$, $U$ and $A$ defined above.

\end{thm}

Following this result, we can define the colored Jones polynomial of any KTG once we fix the value of any colored $\Theta$ graph and describe how it changes under the the operations described above.
\begin{dfn}

~\cite{LV} The colored Jones polynomial of a KTG $\Gamma $ with coloring $\sigma$, denoted by $\langle\Gamma,\sigma\rangle$, is defined by the four equations below.

\[
\langle\Theta; a, b, c\rangle= O^{\frac{a+b+c}{2}}\left[
\begin{matrix}
 &\frac{a+b+c}{2}& &\\
\frac{-a+b+c}{2}&\frac{a-b+c}{2}&\frac{a+b-c}{2}&
\end{matrix}
\right]
\]

\[
\langle F_{\pm}^{e}(\Gamma), \sigma\rangle= f(\sigma(e))^{\pm 1}\langle\Gamma, \sigma\rangle
\]

\[
\langle U^{e}(\Gamma), \sigma\rangle= \langle \Gamma, \sigma\rangle \sum_{\sigma(e)}\frac{O^{\sigma(e)}}{\langle\Theta; \sigma(e), \sigma(b), \sigma(d)\rangle}
\]

\[
\langle A^{\omega}(\Gamma), \sigma\rangle= \langle \Gamma, \sigma\rangle \Delta(a, b, c, \alpha, \beta, \gamma)
\]

\end{dfn}
Particularly, a knot is considered as a $0$-frame KTG without vertices, so we define its colored Jones polynomial as $J_{K}(n+1)=(-1)^{n}\langle K, n \rangle$, where $n$ denotes the color of the single edge of the knot, and the $(-1)^n$ term is to normalize the unknot as $J_{O}(n)=[n]$

In above formulas, the quantum integer $[k]=\frac{v^{2k}-v^{-2k}}{v^2-v^{-2}}$, $[k]!=[k][k-1]\ldots[1]$, the symmetric multinomial coefficieat is defined as:
\[\left[
\begin{matrix}
 a_1+ a_2+ \ldots a_r\\
 a_1,a_2,\ldots, a_r
\end{matrix}
\right]=\frac{[a_{1}+a_{2}+...+a_{r}]!}{[a_{1}]!\ldots[a_{r}]!}.
\]
The value of the $k$-colored unknot is $O^k=(-1)^{k}[k+1]=\langle O,k \rangle $.
The $f$ of the framing change is defined as:
$f(a)=(\sqrt{-1})^{-a} v^{-\frac{1}{2}a(a+2)}$.
The summation in the equation of unzip is taken over all admissible colorings of the edge $e$ that has been unzipped.
$\Delta$ is the quotient of the $6j$-symbol and the $\Theta$, \[\Delta(a,b,c,\alpha,\beta,\gamma)=\Sigma\frac{(-1)^z}{(-1)^{\frac{a+b+c}{2}}}\left[
\begin{matrix}
z+1 \\
 \frac{a+b+c}{2}
\end{matrix}
\right]\left[\begin{matrix}
 \frac{-a+b+c}{2}\\
z-\frac{a+\beta+\gamma}{2}
\end{matrix}
\right]\left[\begin{matrix}
 \frac{a-b+c}{2}\\
z-\frac{\alpha+b+\gamma}{2}
\end{matrix}
\right]\left[\begin{matrix}
 \frac{a+b-c}{2}\\
z-\frac{\alpha+\beta+c}{2}
\end{matrix}
\right].
\]
The summation range of $\Delta$ is indicated by the binomials. Note that this $\Delta$ is not the one in Theorem 2.4, we just maintain their traditional notations and this won't cause any ambiguity according to different contexts.

The above definition agrees with the integer normalization used in~\cite{Cos}, which shows that $\langle \Gamma, \sigma\rangle$ is a Laurent polynomial in $v$ and is independent of the choice of operations to produce the KTG.

As illustrated in Figure 3, we obtain the colored Jones polynomial of the knot $K=M(\frac{1}{r},\frac{1}{s-\frac{1}{u}},\frac{1}{t})$ as follows. Starting from a $\theta$ , we first apply three $A$ moves, then four $F$ moves on the edges $a$, $b$, $c$, $d$ and then unzip these four twisted edges. Finally, to get the $0$-frame colored Jones polynomial we need to multiply the factor $f(n)^{-2(r+s+u+t)-2\cdot writhe}$ (in this case $f(n)^{-4u}$) to cancel the framing produced by the operations and the writhe of the knot.

\begin{figure}[!ht]
\centering

\includegraphics{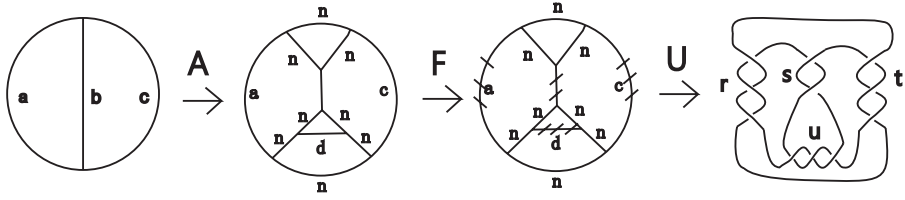}

\caption{The operations to produce the knot $K$ from a $\theta$. For simplicity, in this example we set $r=u=-3, s=2, t=3$. }
\end{figure}

\begin{lem}

The colored Jones polynomial of the Montesinos knot $K=M(\frac{1}{r},\frac{1}{s-\frac{1}{u}},\frac{1}{t} ) $ is
\[
\begin{split}
  J_{K}(n+1)&=(-1)^{n}f^{-4u}(n)\sum_{a, b, c, d\in D_n}^{}
  \langle\Theta;a, b, c\rangle\Delta^{2}(a, b, c, n, n, n)\Delta(b, n, n, d, n, n)f^{r}(a)f^{s}(b)\\&f^{t}(c)f^{u}(d)O^{a}O^{b}O^{c}O^{d}
  \langle\Theta;a, n, n\rangle^{-1}\langle\Theta;b,n,n\rangle^{-1}\langle\Theta;c, n, n\rangle^{-1}\langle\Theta;d, n, n\rangle^{-1},
\end{split}
\]
where the domain $D_n$ is defined such that $a$,$b$,$c$,$d$ are all even with $0\leq a, b, c,d \leq 2n$, and $a$,$b$,$c$ satisfy the triangle inequality.

\end{lem}

To calculate the degree of the colored Jones polynomial,  we need to analyze the the factors of the summands. The following lemma is from~\cite{LV}.

\begin{lem}~\cite{LV}

\[
d_{+}\langle\Theta; a, b, c\rangle= a(1-a)+b(1-b)+c(1-c)+\frac{(a+b+c)^{2}}{2},
\]

\[
d_{+}\langle F_{\pm}^{e}(\Gamma), \sigma\rangle= \pm d_{+}f(\sigma(e))\langle\Gamma, \sigma\rangle,
\]

\[
d_{+}\langle U^{e}(\Gamma), \sigma\rangle\geq d_{+}\langle\Gamma, \sigma\rangle+max_{\sigma(e)}(d_{+}O^{\sigma(e)}-d_{+}\langle\Theta; \sigma(e), \sigma(b), \sigma(d)\rangle),
\]

\[
d_{+}\langle A^{\omega}(\Gamma), \sigma\rangle= d_{+}\langle \Gamma, \sigma\rangle+ d_{+}\Delta(a, b, c, \alpha, \beta, \gamma).
\]

\[
\begin{split}
d_{+}\Delta(a, b, c, \alpha, \beta, \gamma)&= g(m+1,\frac{a+b+c}{2}+1)+ g(\frac{-a+b+c}{2},m-\frac{a+\beta+\gamma}{2})\\
&+g(\frac{a-b+c}{2},m-\frac{\alpha+b+\gamma}{2})+g(\frac{a+b-c}{2},m-\frac{\alpha+\beta+c}{2}),\\
\end{split}
\]
where $g(n,k)=2k(n-k)$ and $2m=a+b+c+\alpha+\beta+\gamma-max(a+\alpha,b+\beta,c+\gamma)$

\end{lem}

Now we are ready to prove Theorem 2.4.

\begin{proof}(of Theorem 2.4)
The maximal degree of $J_{K}(n+1)$ satisfies the inequality below.
\[d_{+}J_{K}(n+1)\leq max_{a, b, c, d\in D_{n}}^{}\Phi(a, b, c, d),\]
where
$\Phi(a, b, c, d)= d_{+}\langle\Theta; a, b, c\rangle+ 2d_{+}\Delta(a, b, c, n, n, n)+d_{+}\Delta(b, n, n, d, n, n)+rd_{+}f(a)+sd_{+}f(b)+td_{+}f(s)+ud_{+}f(d)
+d_{+}O^{a}+d_{+}O^{b}+d_{+}O^{c}+d_{+}O^{d}-d_{+}\langle\Theta; a, n, n\rangle-d_{+}\langle\Theta; b, n, n\rangle-d_{+}\langle\Theta; c,n,n\rangle-d_{+}\langle\Theta; d, n, n\rangle-4ud_{+}f(n)$.

Generally, solving $max_{a,b,c,d\in D_{n}}\Phi(a,b,c,d)$ is a problem of quadratic integer programming, which is quite a involved topic~\cite{GvdV}. In this case however, it can be solved by observing its monotonicity.

Note that the feasible region $D_{n}$ is an even integer lattice in a convex polytope in  $ \mathbb{R}^4$ and can be divided into 6 subregions corresponding to 6 different forms of $\Phi$, that is, $D_{n}=D_{n}^{a,b+d}\cup D_{n}^{b,b+d}\cup D_{n}^{c,b+d}\cup D_{n}^{a,2n}\cup D_{n}^{b,2n}\cup D_{n}^{c,2n}$, where $D_{n}^{a,b+d}$, $D_{n}^{b,b+d}$, $D_{n}^{c,b+d}$, $D_{n}^{a,2n}$, $D_{n}^{b,2n}$ and $D_{n}^{c,2n}$ are defined to be $D_{n}\cap \{(a,b,c,d)\mid a\geq b,c;  \ b+d\geq 2n\}$, $D_{n}\cap \{(a,b,c,d)\mid b\geq a,c;\ b+d\geq 2n\}$, $D_{n}\cap \{(a,b,c,d)\mid c\geq b,a;\ b+d\geq 2n\}$, $D_{n}\cap \{(a,b,c,d)\mid a\geq b,c; \ 2n\geq b+d\}$, $D_{n}\cap \{(a,b,c,d)\mid b\geq a,c;\ 2n\geq b+d\}$ and $D_{n}\cap \{(a,b,c,d)\mid c\geq b,a;\ 2n\geq b+d\}$, respectively. Then we calculate the partial derivatives of the real function $\Phi$, and find that on each of the 6 regions we have $\partial_ {d}\Phi >0$, $\partial_{a}\Phi >0$, and in $D_{n}^{a,b+d}$ we have $\partial_{b}\Phi<0 $. So any maximum of $\Phi$ on $D_{n}$ must occur when $d=2n$,  $a=b+c$. Note that in this paper it is more convenient to calculate the real partial derivatives of these even integer functions when we consider their monotonicity, for example, we use $\partial_{d}\Phi >0$ rather than $\Phi(a,b,c,d+2)>\Phi(a,b,c,d)$, obviously the later is implied by the former.

Now we focus on the following 2-variable function $R(b,c)$ in the domain $T_{n}=\{(b,c)\mid b,c\geq 0, b+c\leq 2n\}$.
\begin{equation}
\begin{split}
R(b, c)&=\Phi(b+c, b, c, 2n)\\
&=-\frac{r+s+1}{2}b^{2}-(r+s-1)b-(r+1)bc-\frac{r+t}{2}c^{2}-(r+t-2)c+2un.\\
\end{split}
\end{equation}

Set $A=-\frac{r+s+1}{2}, B=-(r+1), C=-\frac{r+t}{2}, \Delta=4AC-B^2 $. \\

(1)If $ A\geq0 $ or $ C\geq0 $, we have $\partial_{b}R>0$, or $\partial_{c}R>0$. Then the maxima must be on the line $b+c=2n$. Set $Q(b)=R(b,2n-b)$, we have

\[
\begin{split}
\ Q(b)&=R(b,2n-b)\\
&=-\frac{s+t-1}{2}b^{2}+[2(t-1)n-s+t-1]b-2(r+t)n^2-2(r-u+t-2)n.
\end{split}
\]

$Q(b)$ is a quadratic function in $b$ with negative leading coefficient, and its real maximum is at $b_{m}=\frac{2(t-1)n-s+t-1}{s+t-1}$, $b_{m}\in (0,2n)$ for $n$ sufficiently large.
Set $n+1=N=p(\frac{s+t-1}{2})+j$, $0\leq j< \frac{s+t-1}{2}$, then $b_{m}=(t-1)p-1+\frac{2(t-1)j}{s+t-1}$.
Let $b_{0}$ be the even number nearest to $b_{m}$, then we have $b_{0}=(t-1)p-1+v_{j}$, where $v_{j}$ is the odd number nearest to $\frac{2(t-1)j}{s+t-1}$,
 so $b_{0}=\frac{2(t-1)}{s+t-1}N-\frac{2(t-1)}{s+t-1}j+v_{j}-1$. Then we have
\[max_{a,b,c,d\in D_{n}}\Phi(a,b,c,d)= Q(b_{0})=[\frac{2(t-1)^{2}}{s+t-1}-2(r+t)]N^{2}+2(r+u+3)N+c_{j},\]
where $c_j=-\frac{s+t-1}{2}\beta^{2}_{j}-(s+t-1)\beta_j-2(u+2)$, $\beta_j=v_j-1-\frac{2(t-1)}{s+t-1}j$.

Finally, when $b_{m}$ is not an odd integer, the maximum is unique. Otherwise, $\Phi$ has exactly 2 maxima, we have to consider the possibility that the coefficients of the 2 maximal degree terms may cancel out. It is easy to see that for the leading coefficients of the terms of the summation, the $O$ terms contribute $(-1)^{a+b+c+d}=1$, the $f$ terms contribute $(-1)^{\frac{1}{2} (ar+bs+ct+du)}$, $\Theta$ terms contribute $(-1)^{a+b+c+4n+\frac{1}{2}d}=(-1)^{\frac{1}{2}d}$, $\Delta$ terms contribute$(-1)^{\frac{1}{2}b+n}$. Multiplied by the $(-1)^{n}$ in front of the summation, the leading coefficients are $(-1)^{\frac{1}{2}[ra+(s+1)b+tc+(u+1)d]}$. Note that $r,u,t$ are odd integers and $s$ is an even integer, and any maximum must occur on $a=b+c=2n$. So if there are two maximal degree terms, their coefficient are both $1$, and no cancellation will happen.

(2)When $A<0$ and $C<0$, for any fixed $c_0$, $R(b,c_{0})$ is a quadratic function in $b$ with negative leading coefficient, whose axis of symmetry (in the plane $c=c_0$ of the uvR-coordinates) intersects the line $\partial_{b}R=0$ and is perpendicular to the $bc$-plane, so we consider the real value of $R(b,c)$ on the line $\partial_{b}R=0$:
\begin{equation*}
  \left\{
   \begin{array}{l}
  R(b,c)=-\frac{r+s+1}{2}b^{2}-(r+s+1)b-(r+1)bc-\frac{r+t}{2}c^2-(r+t-2)c+2un  \\
  \partial_{b}R=-(r+s+1)b-(r+s-1)-(r+1)c =0.\\
   \end{array}
   \right.
  \end{equation*}
Then we have
\begin{equation}
 \begin{split}
&R\mid _{\partial _b R=0}= R(b,-\frac{r+s+1}{r+1}b-\frac{r+s-1}{r+1})\\
 &=[-\frac{r+s+1}{2(r+1)^{2}}\Delta]b^{2}+\frac{r+s+1}{(r+1)^{2}}[(r+1)(r+t-2)-(r+t)(r+s-1)]b+ const.
 \end{split}
\end{equation}

If $\Delta \neq 0$, $R\mid _{\partial _b R=0}$ is a quadratic function in $b$ whose axis of symmetry is perpendicular to the $bc$-plane at the point
\[P=(b_p, c_p)=(\frac{(r\!+\!1)(r\!+t\!-\!2)\!-\!(r\!+\!t)(r\!+\!s\!-\!1)}{\Delta},\frac{(r\!+\!1)(r\!+\!s\!-\!1)\!-\!(r\!+\!s\!+\!1)(\!r+\!t\!-\!2)}{\Delta}).\]($P$ is actually the intersection of $\partial_{b}R=0$ and $\partial_{c}R=0$).

\begin{figure}[!ht]
\centering

\includegraphics{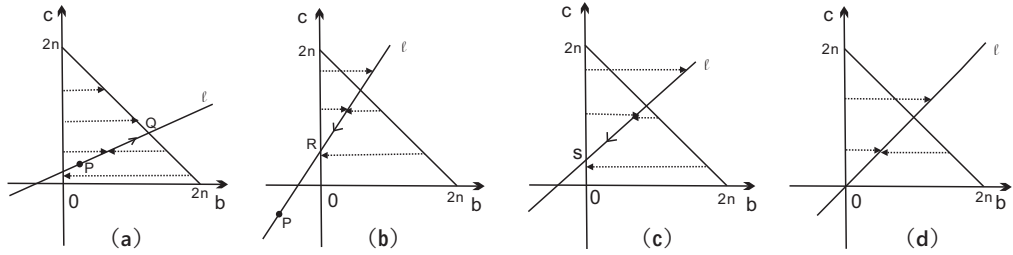}

\caption{$R(b,c)$ is restricted to the triangle domain $T_n$: $(0,0)$-$(0,2n)$-$(2n,0)$, and the arrows indicate its increasing direction; $\ell$ denotes the line $\partial_{b}R=0$.  }
\end{figure}

(2.1)If $A$ and $C<0, \Delta<0$, by Equation (3.2), $R\mid _{\partial _b R=0}$ is a quadratic function in $b$ with positive leading coefficient. And we have $b_p\leq 0$, $c_p\leq 0$. See Figure 4(a). The arrows indicate the increasing direction. For sufficiently large $n$, any maximum must be on the segment $[Q,(0,2n)]$, then the argument will be the same as that of case (1).

(2.2)If $A,C<0$, and $\Delta>0$, $R\mid _{\partial _b R=0}$ is a quadratic function in $b$ with negative leading coefficient. And we have $b_p\geq 0$, $c_p\geq 0$. Any maximum must occur on $OR $. See Figure 4(b). Note that
 $R(0,c)=-\frac{r+t}{2}c^{2}-(r+t-2)c+2un$ decreases in $[0,+\infty)$, the maximum must occur on $O=(0,0)$, so
\[d_{+}J_{K}(n+1)=R(0,0)=2un, \ d_{+}J_{K}(n)=2u(n-1).\]

(2.3)If $A,C<0$ ,$\Delta = 0$, and $(r+s-1)^2+(r+t-2)^2\neq 0$, $R\mid _{\partial _b R=0}$ is a decreasing linear function in $b$, any maximum must occur on OS. See Figure 4(c). $R(0,c)$ decreases in $[0,+\infty)$, the maximum must be on $O=(0,0)$, so
\[ d_{+}J_{K}(n)=2u(n-1).\]

(2.4)If $A,C<0$, $\Delta =0$, and $(r+s-1)^{2}+(r+t-2)^{2}=0$, then we immediately have $r=-3, s=4, t=5$, $R(b,c)=-(b-c)^{2}+2un$, the maxima are $R(0,0)=R(2,2)=\ldots =R(k,k)$, where $k=n$ when $n$ is even, and $k=n-1$ when $n$ is odd. See Figure 4(d). By a similar argument with the last paragraph of (1), we conclude that there are no cancellations between the the highest-degree coefficients, so
 \[d_{+}J_{K}(n)=2u(n-1).\]

\end{proof}

\section{Boundary Slope and Euler Characteristic  }

 The main idea of Hatcher-Oertel edgepath system, based on the work of~\cite{HT85}, is to deal with properly embedded surfaces in Montesinos knot complement  combinatorially. For details see~\cite{HO89, IM07}. Briefly speaking, the so called \textit{candidate surfaces} are associated to \textit{admissible edgepath systems} in a diagram $\mathcal{D}$ in the $uv$-plane. The vertices of $\mathcal{D}$ correspond to projective curve systems $[a, b, c]$ on the 4-punctured sphere carried by the train track in Figure 5(a) via $u=\frac{b}{a+b}$, $v=\frac{c}{a+b}$.
\begin{figure}[!ht]
\centering

\includegraphics{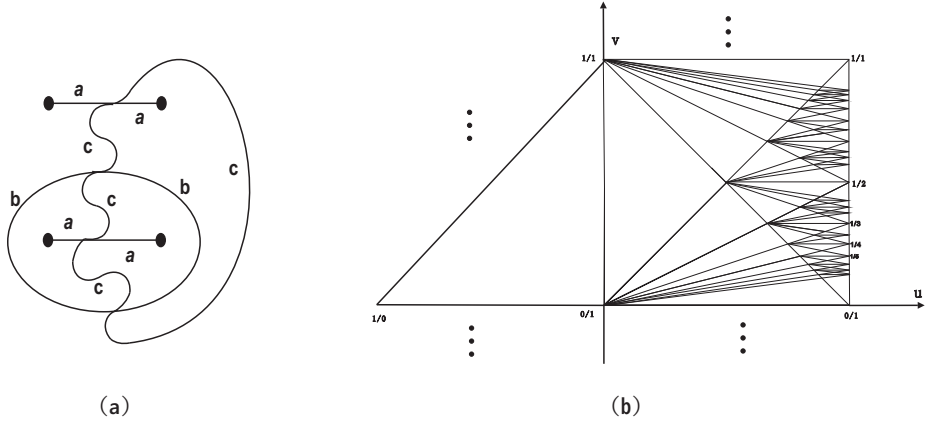}

\caption{(a)The train track in a 4-punctured sphere. (b)The diagram $\mathcal{D}$ in the $uv$-plane. }
\end{figure}

There are three types of  vertices of $\mathcal{D}$:

   (1)the arcs with slope $\frac{p}{q}$ denoted by $\langle\frac{p}{q}\rangle$, corresponding to the projective curve systems $[1, q-1, p]$, with the \textit{uv}-coordinates $(\frac{q-1}{q}, \frac{p}{q})$,

   (2)the circles with slope $\frac{p}{q}$ denoted by $\langle\frac{p}{q}\rangle ^{\circ}$, corresponding to the projective curve systems $[0, p, q]$, with the \textit{uv}-coordinates $(1, \frac{p}{q})$,

   (3)the arcs with slope $\infty$ denoted by $\langle \infty \rangle$, with the \textit{uv}-coordinates $(-1,0)$.

   And there are 6 types of  edges in $\mathcal{D}$:

   (1)the \textit{non-horizontal } edges connecting the vertex $\frac{p}{q}$ to the vertex $\frac{r}{s}$ with $\mid ps- qr\mid =1$, denoted by $[\langle \frac{r}{s}\rangle , \langle \frac{p}{q}\rangle ]$,

   (2)the \textit{horizontal } edges connecting $\langle\frac{p}{q}\rangle ^{\circ}$ to $\langle\frac{p}{q}\rangle $, denoted by $[ \langle \frac{p}{q}\rangle ,\langle\frac{p}{q}\rangle ^{\circ} ]$,

   (3)the \textit{vertical } edges connecting $\langle z\rangle$ to $\langle z\pm 1\rangle$ , denoted by $[z, z\pm 1]$, here $z\in $ $\mathbb{Z}$,

   (4)the \textit{infinity} edges connecting $\langle z\rangle$ to $\langle \infty \rangle$ denoted by $[\infty, z]$,

   (5)the \textit{constant} edges which are  points on the horizontal edge $[ \langle \frac{p}{q}\rangle ,\langle\frac{p}{q}\rangle ^{\circ} ]$ with the form $ \frac{k}{m} \langle \frac{p}{q}\rangle + \frac{m-k}{m}\langle\frac{p}{q}\rangle ^{\circ} $,

   (6)the \textit{partial} edges which are parts of non-horizontal edges $[\langle \frac{r}{s}\rangle , \langle \frac{p}{q}\rangle ]$ with the form $[\frac{k}{m}\langle \frac{r}{s}\rangle + \frac{m-k}{m} \langle \frac{p}{q}\rangle, \langle \frac{p}{q}\rangle ] $.

   An \textit{edgepath} denoted by $\gamma$ in $\mathcal{D}$ is a piecewise linear path beginning and ending at rational points of $\mathcal{D}$. An admissible edgepath system denoted by $\Gamma =(\gamma_{1},\gamma_{2},\ldots, \gamma_{n})$ is an \textit{n}-tuple of edgepaths in $\mathcal{D}$ with following properties.\\
(E1)The starting point of $\gamma_i$ is on the horizontal edge $[\langle \frac{p_i}{q_i} \rangle, \langle \frac{p_i}{q_i} \rangle^{\circ}]$, and if it is not the vertex $\langle \frac{p_i}{q_i}\rangle$, $\gamma_i$ is constant.\\
(E2)$\gamma_i$ is minimal, that is, it never stops or retraces itself, nor does it ever go along two sides of the same triangle of $\mathcal{D}$ in succession.\\
(E3)The ending points of $\gamma_i$'s are rational points of $\mathcal{D}$ with their $\textit{u}$-coordinates  equal and $\textit{v}$-coordinates adding up to zero.\\
(E4)$\gamma_i$ proceeds monotonically from right to left, ``monotonically'' in the weak sense that motion along vertical edges is permitted.

In~\cite{HO89}, a finite number of candidate surfaces are associated to the each admissible edgepath system, then every essential surface in knot complement with non-empty boundary of finite slope is isotopic to one of the candidate surfaces. Finally, to rule out the inessential surfaces, the notion of \textit{r}-value is developed in ~\cite{HO89}, in our case however, we only need the following convenient criterion from~\cite{Ich}.

\begin{lem}~\cite{Ich}
   For a admissible edgepath system having ending points with positive $u$-coordinate, if all the last edges of the edgepaths travel in the same direction (upward or downward) from right to left, then all the candidate surfaces associated to the edgepath system are essential.
   \end{lem}

The boundary slope of a essential surface $S$ is computed by $\tau(S)- \tau(S_0)$, where $\tau(S)$ is the \textit{total number of twist} (or twist for short) of an essential surface $S$, and $S_0$ is the Seifert surface. For a candidate surface $S$ associated to the admissible edgepath system $\Gamma$, we have~\cite{IM07}
\begin{equation}
  \tau(S)=\sum_{\gamma_{i}\in \Gamma_{non-const}} \sum_{e_{i,j}\in \gamma_{i}}-2\sigma(e_{i,j})|e_{i,j}|.
\end{equation}
 In above formula, $|e|$ is the \textit{length} of an edge $e$, which is defined to be $0$, $1$, or $\frac{k}{m}$ for a constant edge, a complete edge or a partial edge $[\frac{k}{m}\langle \frac{r}{s}\rangle + \frac{m-k}{m} \langle \frac{p}{q}\rangle, \langle \frac{p}{q}\rangle ] $, respectively. And $\sigma(e)$ is the \textit{sign} of a non-constant edge $e$, which is defined to be $+1$ or $-1$ according to whether the edge is increasing or decreasing (from right to left in \textit{uv}-plane) respectively for a non-$\infty$ edge; for an $\infty$ edge the sign is defined to be $0$.

Now we are ready to prove Theorem 2.5.

\begin{proof}(of Theorem 2.5)
By the method of~\cite{HO89} (pg461), when $r,u,t$ are odd and $s$ is even with $u\leq-1$, $r<-1<1<s,t$, we directly find the edgepath system of the Seifert surface $S_0$:

\[
\begin{split}
&\delta_{1}:[\langle0\rangle,\langle\frac{1}{r}\rangle],\\
&\delta_{2}:[\langle0\rangle, \langle\frac{1}{s+1}\rangle, \ldots,\langle\frac{-u-i}{s(-u-i)+1}\rangle, \ldots, \langle\frac{-u}{-su+1}\rangle],\\
&\delta_{3}:[\langle0\rangle, \langle\frac{1}{t}\rangle],\\
\end{split}
\]
where $0\leq i \leq u-1$.
$S_0$ has boundary slope 0, and by the formula (3.4) from~\cite{IM07},
\[-\frac{\chi(S_0)}{\sharp S_0}=2-u-2,\]
\[\frac{\chi(S_0)}{\sharp S_0}=u.\]
By now we have proved (2).

For (1), with direct calculations we always have $\Delta<0$, and we claim that there exists an admissible edgepath system having ending points with $\textit{u}$-coordinate  $u_{0}=\frac{(t-1)s}{ts+t-1}$ in $\textit{uv}$-plane. In fact, $u_0$ is just the solution of the equation $v_{1}(u)+v_{2}(u)+ v_{3}(u)=0$, where the linear functions $v=v_{1}(u)$, $v_{2}(u)$ and $v_{3}(u)$ are determined by the lines through the edges $[\langle0\rangle,\langle\frac{1}{t}\rangle]$, $[\langle0\rangle,\langle\frac{1}{s+1}\rangle]$ and $[\langle -1\rangle,\langle -\frac{1}{2}\rangle...,\langle\frac{1}{r}\rangle]$, respectively. Denote by $u_t$, $u_{s+1}$ and $u_r$ the $\textit{u}$-coordinates of $\langle \frac{1}{t}\rangle$, $\langle \frac{1}{s+1}\rangle$ and $\langle \frac{1}{r}\rangle$ respectively. With direct calculations we have
\[u_0-u_t=-\frac{(t-1)^2}{t(st+t-1)}<0,\]
 \[u_0-u_{s+1}=-\frac{s^2}{(s+1)(st+t-1)}<0,\]
  \[u_0-u_r=-\frac{\Delta}{r(st+t-1)}<0,\]
  so $u_o$ must be on the left of $u_t$, $u_{s+1}$ and $u_r$. Suppose the edgepath of the $\frac{1}{r}$-tangle ends on edge $[\frac{1}{r+k+1}, \frac{1}{r+k}]$, where $0\leq k\leq -r-2$, then $u_0$ must be the $\textit{u}$-coordinate of ending points of the admissible edgepath system $\Gamma$ below:
\[
\begin{split}
&\gamma_{1}:[(\frac{(t\!-\!1\!)^{2}}{s+t-\!1}\!-\!r\!-\!t-\!k)\langle\frac{1}{r\!+\!k\!+\!1\!}\!\rangle\!+\!(\!1\!+\!r\!+\!t\!+\!k-\!\frac{(\!t\!-\!1)^{2}}{s\!+\!t\!-\!1})\!\langle\frac{1}{r+k}\rangle,
 \langle\frac{1}{r\!+\!k}\!\rangle, \!\langle\frac{1}{r\!+\!k\!-1}\!\rangle, \ldots, \!\langle\frac{1}{r}\!\rangle], \\
&\gamma_{2}:[\frac{s}{s+t-1}\langle0\rangle+\frac{t-1}{s+t-1}\langle\frac{1}{s+1}\rangle, \langle\frac{1}{s+1}\rangle, \ldots, \langle\frac{-u-i}{s(-u-i)+1}\rangle,
  \ldots, \langle\frac{-u}{-su+1}\rangle],\\
&\gamma_{3}:[\frac{t-1}{s+t-1}\langle0\rangle+\frac{s}{s+t-1}\langle\frac{1}{t}\rangle, \langle\frac{1}{t}\rangle],\\
\end{split}
\]
where $0\leq i \leq -u-1$, and the length of the partial edges are calculated via $u_0$ by formula (3.1) from~\cite{IM07}. By Lemma 4.1, any candidate surface associated to $\Gamma$ must be essential.

By formula (4.1), the twist of the surface $S$ associated to the above edgepath system is
\[
\begin{split}
\tau(S)&=\sum_{r_{i}\in \Gamma_{non-const}^{}} \sum_{e_{i,j}\in r_{i}}-2\sigma(e_{i,j})|e_{i,j}|
\\&=2[\frac{(t-1)^{2}}{s+t-1}-r-t-k]+2k+\frac{2s}{s+t-1}+2(-u-1)+\frac{2(t-1)}{s+t-1}\\
&=\frac{2(t-1)^{2}}{s+t-1}-2(u+r+t).
\end{split}
\]
The twist of the Seifert surface $S_0$ is
\[\tau (S_0)=-2-2u+2=-2u.\]
So the boundary slope of $S$ is
\[\tau(S)-\tau(S_0)=\frac{2(t-1)^{2}}{s+t-1}-2(r+t).\]
Finally, by the formula (3.5) from~\cite{IM07} we have
\[
\begin{split}
 -\frac{\chi(S)}{\sharp S}&=\sum_{r_{i}\in \Gamma_{non-const}}|r_{i}|+N_{const}-N+(N-2-\sum_{r_{i}\in \Gamma_{const}}\frac{1}{q_{i}})\frac{1}{1-u_{0}}\\
&=\sum_{r_{i}\in\Gamma_{non-const}}|r_{i}|-3+\frac{1}{1-u_{0}}\\
&=\frac{(t-1)^{2}}{s+t-1}-(u+r+t)-3+\frac{ts+t-1}{s+t-1}\\
&=-(u+r+3).\\
\end{split}
\]

\[\frac{\chi(S)}{\sharp S}=u+r+3.\]

\end{proof}

\bibliographystyle{amsplain}

\end{document}